\documentclass[a4paper,11pt]{article}
\usepackage[utf8]{inputenc}

\usepackage{cite}
\usepackage{graphicx}
\usepackage{graphicx}
\usepackage{pstricks}
\usepackage{appendix}
\usepackage{dsfont}
\usepackage{amsmath}
\usepackage{amssymb}
\usepackage{amsthm}
\usepackage{mathtools}
\usepackage{comment}
\usepackage[shortlabels]{enumitem}
\usepackage{tikz,lipsum,lmodern}
\usepackage[most]{tcolorbox}

\usepackage{xcolor}
\usepackage{blindtext}

\usepackage[nottoc]{tocbibind}

\usepackage[linktoc=page, colorlinks, linkcolor=blue, citecolor=blue]{hyperref}

\numberwithin{equation}{section}
\newtheorem{theorem}{Theorem}[section]
\newtheorem{corollary}[theorem]{Corollary}

\newtheorem{lemma}[theorem]{Lemma}

\theoremstyle{remark}
\newtheorem{remark}[theorem]{Remark}

\DeclareMathOperator{\tr}{Tr}

\DeclareMathOperator{\rad}{rad}

\newcommand{\R}{\mathbf R}

\newcommand{\E}{{\mathbb E}}
\newcommand{\EE}{{\mathcal E}}
\newcommand{\FF}{{\mathcal F}}

\DeclarePairedDelimiter \abs{\lvert}{\rvert} 
\DeclarePairedDelimiter \norm{\lVert}{\rVert}
\DeclarePairedDelimiterX \ip[2]{\langle}{\rangle}{#1,#2}
\DeclarePairedDelimiterXPP \Prob[1]{\mathbb{P}}\{\}{}{ #1} 
\DeclarePairedDelimiterXPP \Probevent[1]{\mathbb{P}}(){}{#1} 

\def \R {\mathbb{R}}
\def \tran {\mathsf{T}}
\def \psitwo {{\psi_2}}
\def \psione {{\psi_1}}
\def \one {{\textbf 1}}

\title{On the Dimension-Free Concentration of Simple Tensors via Matrix Deviation}

\begin{document}
\author{
  Pedro Abdalla
  \qquad
  Roman Vershynin\thanks{Department of Mathematics, UC Irvine}
}

\maketitle

\begin{abstract}
    We provide a simpler proof of a sharp concentration inequality for subgaussian simple tensors obtained recently by Al-Ghattas, Chen, and Sanz-Alonso. Our approach uses a matrix deviation inequality for $\ell^p$ norms and a basic chaining argument.
\end{abstract}

\section{Introduction}
Let $X$ be a mean-zero random vector in $\R^d$ with covariance matrix $\Sigma \coloneqq \E XX^\tran \in \mathbb{R}^{d\times d} $, and let $p \ge 2$ be an integer. We define the simple tensor $X^{\otimes p}$ as the $p$-order tensor whose entries are the products of all $p$-tuples of the entries of $X$. Can we estimate the expected value of the simple random tensor $X^{\otimes p}$ from a sample $X_1,\ldots,X_N$ of i.i.d. copies of $X$? The simplest estimator is the empirical mean. The error of this estimator is
\begin{equation}    \label{eq:main_problem}
    \norm*{\frac{1}{N}\sum_{i=1}^N X_i^{\otimes p}- \E X^{\otimes p}}
    \coloneqq \sup_{v\in S^{d-1}} \abs*{\frac{1}{N}\sum_{i=1}^N \ip{X_i}{v}^p - \E\ip{X}{v}^p}.
\end{equation}

How large is it? We are interested in dimension-free bounds -- those that do not depend on the full dimension $d$, but instead on the ``effective dimension'' of the distribution, captured by the effective rank of the covariance matrix:
\begin{equation*}
    r(\Sigma)
    \coloneqq \frac{\tr(\Sigma)}{\norm{\Sigma}}.
\end{equation*}
This question is nontrivial already for $p=2$, where it becomes the covariance estimation problem. The standard $\varepsilon$-net argument gives optimal dimension-dependent bounds (see \cite[Chapter 4]{vershynin2018high}). But getting dimension-free bounds even for $p=2$ is harder. They were first obtained by Lounici and Koltchinskii \cite{koltchinskii2017concentration} for subgaussian distributions using generic chaining. A simpler argument for Gaussian distribution was found by Van Handel \cite{van2017structured}, while Liaw, Mehrabian, Plan and the second author \cite{liaw2017simple} show how to deduce the result for all subgaussian distributions from their matrix deviation inequality (see \cite[Chapter 9]{vershynin2018high} for an introduction). That last approach also lets one take the supremum in \eqref{eq:main_problem} over any bounded set $T \subset \R^d$ -- the setting we consider in the present paper, too.

Much less has been known when $p>2$. Optimal bounds that depend on dimension were proved in \cite{guedon2007lp,vershynin_lp_marginals,adamczak2010quantitative} when $p$ is even. Very recently, Al-Ghattas, Chen and Sanz-Alonso \cite{al2025sharp} proved optimal bounds for $p\ge 2$. Their proof relies on a local chaining argument due to Mendelson \cite{mendelson2016upper} combined with an intricate analysis of coordinate projections due to Bednorz \cite{bednorz2014concentration}. 

The goal of our work is to provide a simpler proof of this fact by exploiting tools from non-asymptotic random matrix theory. Specifically, we can leverage the matrix deviation inequality (see Theorem \ref{thm:matrix_dev}) to bypass Bendorz's analysis of coordinate projections, and we can run the chaining argument for the simpler $L^2$ structure of the process and bypass the delicate chaining argument of Mendelson.

To state the result, we say that a random vector in $\mathbb{R}^d$ is isotropic if its covariance matrix is the identity matrix $I_d$. Next, let $g\sim N(0,I_d)$ be the standard multivariate Gaussian random vector. For a set $T \subset \R^d$, the radius and the Gaussian complexity are defined by
\begin{equation*}
    \rad(T) \coloneqq \sup_{v\in T}\norm{v}_2, \quad \gamma(T) \coloneqq \E \sup_{v\in T} \abs{\ip{g}{v}}.
\end{equation*}
Notice that $\sqrt{2/\pi}\rad(T) \le \gamma(T)$ by Jensen's inequality.
Here and in the rest of the paper, the $\psi_2$ norm and the $L^p$ norm of a random variable $Z$ are defined by 
$$
\norm{Z}_\psitwo 
\coloneqq \inf \left\{ t>0: \E e^{Z^2/t^2} \le 2 \right\}
\quad \text{and} \quad
\norm{Z}_{L^p} = (\E \abs{Z}^p)^{1/p}.
$$
We say that 
a random vector $X$ is subgaussian if there is an $K>0$ such that
$$
\norm{\ip{X}{v}}_\psitwo 
\le K \norm{\ip{X}{v}}_{L^2},
\quad \text{for any } v\in S^{d-1}.
$$
For basics on subgaussian distributions, see \cite{vershynin2018high}. The following theorem and corollary were proved by Al-Ghattas, Chen and Sanz-Alonso \cite{al2025sharp}. We provide a simpler proof of these facts.

\begin{theorem} \label{thm:main}
    Let $Z_1,\ldots,Z_N$ be i.i.d. copies of a mean-zero subgaussian isotropic random vector $Z$, and let $T\subset \R^d$ be a bounded set. Then, for any integer $p\ge 2$, we have
    $$
    \E\sup_{v\in T}\abs*{\sum_{i=1}^N \ip{Z_i}{v}^p - N \E \ip{Z}{v}^p} 
    \lesssim_{K,p} \gamma(T)^p + \sqrt{N}\gamma(T)\rad(T)^{p-1},
    $$
    where the sign $\lesssim_{K,p}$ hides factors that depend only on $K$ and $p$.
\end{theorem}

\begin{corollary}   \label{cor: anisotropic}
    Let $X$ be a mean-zero subgaussian random vector with covariance matrix $\Sigma$ and $X_1,\ldots,X_N$ be i.i.d. copies of $X$. Then, for any integer $p\ge 2$, we have
    $$
    \E \norm*{\frac{1}{N}\sum_{i=1}^N X_i^{\otimes p}- \E X^{\otimes p}} 
    \lesssim_{K,p} \norm{\Sigma}^{p/2}\left(\frac{r(\Sigma)^{p/2}}{N}+\sqrt{\frac{r(\Sigma)}{N}}\right).
    $$
    
\end{corollary}

To derive Corollary \ref{cor: anisotropic}, write $X=\Sigma^{1/2}Z$ for some isotropic random vector $Z$ and apply Theorem \ref{thm:main} for the ellipsoid $T=\Sigma^{1/2}S^{d-1}$, noting that
$$
\rad(\Sigma^{1/2}S^{d-1}) = \norm{\Sigma}^{1/2},
\quad
\gamma(\Sigma^{1/2}S^{d-1})\leq \tr(\Sigma)^{1/2} = \norm{\Sigma}^{1/2}r(\Sigma)^{1/2}.
$$

\section{Matrix Deviation Approach}
Our approach will also exploit the matrix deviation inequality, as was first done for the case $p=2$ in \cite{liaw2017simple}, see \cite[Chapter 9]{vershynin2018high}. We will use the following version for $p \ge 2$ proved by Sheu and Wang \cite{sheu2023matrix}:

\begin{theorem}[$\ell^p$ matrix deviation inequality]   \label{thm:matrix_dev}
    Let $p\ge 2$. Let $A \in \R^{N\times d}$ be a random matrix whose rows are i.i.d. copies of a mean-zero isotropic subgaussian random vector $Z \in \mathbb{R}^d$. Consider the stochastic process
    \begin{equation*}
        Z_v
        \coloneqq \abs*{\norm{Av}_{p}-N^{1/p}\norm{\ip{Z}{v}}_{L^p}}, \quad v \in \R^d.
    \end{equation*}
Then $Z_v$ has Lipschitz subgaussian increments:
    \begin{equation*}
        \norm{Z_v-Z_u}_\psitwo 
        \lesssim_{K,p} \norm{u-v}_2
        \quad \text{for every } u,v \in \R^d.
    \end{equation*}
\end{theorem}

\begin{remark}[A bound on the process]
\label{rmk:tail_matrixdev}
It immediately follows from Theorem \ref{thm:matrix_dev} combined with Talagrand's majorizing measure theorem (see \cite[Theorem 8.5.5]{vershynin2018high}) that for any bounded set $T \subset \R^n$ and any $u>0$, we have 
\begin{equation*}
\sup_{v \in T} |Z_v| \lesssim_{K,p} \gamma(T) + u\rad(T)
\end{equation*}
with probability at least $1-2e^{-u^2}$.
\end{remark}

\begin{remark}[Subgaussianity]
    Theorem \ref{thm:matrix_dev} is stated in \cite{sheu2023matrix} under a stronger assumption that the entries $A$ are i.i.d. subgaussian, but the proof actually only relies on the rows being subgaussian.
\end{remark}

As we mentioned, for $p=2$ one can derive Theorem \ref{thm:main} directly from the matrix deviation inequality. However, for $p>2$, in addition to matrix deviation, we will need an extra chaining step. 

To run the chaining, we need a couple of standard bounds on the order statistic of a sequence of independent subgaussian random variables. Given real numbers $X_1,\ldots,X_n$, we denote by $X^*_1, \ldots, X^*_n$ a nonincreasing rearrangement of the numbers $\abs{X_1}, \ldots, \abs{X_n}$, so that $\abs{X^*_1} \ge \abs{X^*_2} \ge \cdots \ge \abs{X^*_n}$. We define the function $\ln_+(x)\coloneqq\max\{\ln(x),0\}$.

\begin{lemma}[Order statistics]   \label{lem: rearrangement}
    Let $X_1,\ldots,X_n$ be independent random variables satisfying $\norm{X_i}_\psitwo \le 1$ for all $i$. Let $t>0$, $q \ge 2$ and 
    $k \coloneqq t/\ln_+(en/t)$.
    Then,\footnote{By convention, we set $1/0=\infty$ in defining $k$, and the sum over an empty set is zero.} with probability at least $1-2e^{-t}$, we have
    \begin{equation} \label{eq: order statistics}
        \sum_{i \le 3k} (X_i^*)^2 \lesssim t
        \quad \text{and} \quad 
        \sum_{i > k} (X_i^*)^q \lesssim_q n.
    \end{equation}
\end{lemma}
We postpone the proof to Appendix \ref{a: order statistics} and prove the main result.

\begin{proof}[Proof of Theorem \ref{thm:main}]

{\em Step 1: Symmetrization.}
Denote by $\mathcal{F}$ be the class of functions $f:\mathbb{R}^d\rightarrow \mathbb{R}$ of the form $f(Z)=\langle Z,v\rangle$, where $v\in T$. 
Then, by the classical Gine-Zinn symmetrization (see for example \cite[Chapter 6]{vershynin2018high}), we have
\begin{align}    
    \E\sup_{v \in T} \abs[\bigg]{\sum_{i=1}^N \ip{Z_i}{v}^p - N \E \ip{Z}{v}^p}
    &= \E \sup_{f \in \FF} \abs[\bigg]{\sum_{i=1}^N \big( f^p(Z_i) - \E f^p(Z_i) \big)} \nonumber\\ 
    &\leq 2 \E \sup_{f \in \FF} \abs[\bigg]{\sum_{i=1}^N \varepsilon_i f^p(Z_i)}, \label{eq: symmetrization}
\end{align}
where $\varepsilon_1,\ldots,\varepsilon_N$ are i.i.d. Rademacher random variables, independent of $(Z_i)$ and the expectation is taken with respect to both the law of $(Z_i)$ and $(\varepsilon_i)$. 

\smallskip

{\em Step 2: Generic chaining.}
Given an admissible sequence $(T_s)_{s\ge 0}$ in $T$, that is, a collection of sets $(T_s)_{s\ge 0}\subset T$ satisfying $|T_s|\le 2^{2^s}$ and $|T_0|=1$, we consider the collection of maps $\pi_s:T\rightarrow T_s$ that map $v\in T$ to a nearest point $\pi_s v\in T_s$. Denote by $\mathcal{F}_s$ the class of functions from $\mathbb{R}^d$ to $\mathbb{R}$ defined by $(\pi_s f)(Z) \coloneqq \langle Z,\pi_sv\rangle$, where $v \in T$.

By Talagrand's majorizing measure theorem (see \cite[Chapter 8]{vershynin2018high} or \cite[Chapter 2]{talagrand2022upper}), there exists an admissible sequence such that
\begin{equation} \label{eq:Talagrand_MM}
\sup_{f\in \mathcal{F}}\sum_{s\ge 0}2^{s/2} \norm{\Delta_s f}_{L^2(Z)}=\sup_{v\in T}\sum_{s\ge 0}2^{s/2} \norm{\Delta_s v}_2 \lesssim \gamma(T),
\end{equation}
where $\Delta_s f \coloneqq \pi_{s+1}f-\pi_s f$ and $\Delta_s v \coloneqq \pi_{s+1}v-\pi_s v$.
Fix this admissible sequence. We can decompose any $f\in \mathcal{F}$ along this chain: 
\begin{equation*}
    f^p(Z) = \sum_{s\ge 0} \big( (\pi_{s+1}f)^p-(\pi_sf)^p \big)(Z) + (\pi_{0}f)^p(Z),
\end{equation*}
where the series converges in $L^2(\mathbb{R}^d)$.
Decomposing $f^p(Z_i)$ this way gives
\begin{equation}    \label{eq: complete_term}
    \sum_{i=1}^N \varepsilon_i f^p(Z_i)
    = \sum_{s\ge 0}\sum_{i=1}^N \varepsilon_i \big( (\pi_{s+1}f)^p-(\pi_sf)^p \big)(Z_i) + \sum_{i=1}^N\varepsilon_i(\pi_{0}f)^p(Z_i).
\end{equation}
Now, $\pi_0 f \in \FF_0$ and $\FF_0$ is a singleton (as is $T_0$), so $\pi_0 f$ does not depend on the choice of $f$. Thus, if we take supremum over $f \in \FF$ and then take expectation, we can bound \eqref{eq: symmetrization} as follows:
\begin{align}
    \E &\sup_{f \in \FF} \abs[\bigg]{\sum_{i=1}^N \varepsilon_i f^p(Z_i)} \nonumber\\
    &\lesssim \E \sup_{f \in \FF} \abs[\bigg]{\sum_{s\ge 0}\sum_{i=1}^N \varepsilon_i \big( (\pi_{s+1}f)^p-(\pi_sf)^p \big)(Z_i)} + \E \abs[\bigg]{\sum_{i=1}^N\varepsilon_i(\pi_{0}f)^p(Z_i)}.  \label{eq: chained}
\end{align}

\smallskip

{\em Step 3: Hoeffding's inequality.}
The main work goes into bounding the first term on the right hand side of \eqref{eq: chained}. To do so, let's first fix $s \ge 0$ and $f \in \mathcal{F}$. Split the sum $\sum_{i=1}^N$ into $\sum_{i\in I_s}$ and $\sum_{i\notin I_s}$, where the set  $I_s\subset [N]$ will be chosen later independently of the signs $\varepsilon_i$. To bound $\sum_{i\in I_s}$, use the triangle inequality. To bound $\sum_{i \notin I_s}$, use Hoeffding's inequality conditionally on $(Z_i)$ with the choice $t=u2^{s/2}$ (see for example \cite[Chapter 2]{vershynin2018high}).

For any $u> 0$, we get that with probability at least $1-2e^{-2u^2 2^{s}}$,
\begin{align}
    \abs[\bigg]{&\sum_{i=1}^N \varepsilon_i \big( (\pi_{s+1}f)^p-(\pi_sf)^p \big)(Z_i)} \nonumber\\
    &\lesssim \underbrace{\sum_{i\in I_s} \abs*{(\pi_{s+1}f)^p-(\pi_sf)^p}(Z_i)}_{\eqqcolon A_s(f)} + \underbrace{u2^{s/2}\bigg(\sum_{i\notin I_s} \big( (\pi_{s+1}f)^p-(\pi_sf)^p \big)^2(Z_i)\bigg)^{1/2}}_{\eqqcolon B_s(f)}. \label{eq: As Bs}
\end{align}
By the numeric inequality $|a^p-b^p|\le |a-b|\big(|a|^{p-1}+|b|^{p-1}\big)$, we have
\begin{align*}
    A_s(f) &\lesssim_p \sum_{i\in I_s} \abs*{(\Delta_s f)(Z_i)} \Big( \abs{(\pi_{s+1}f)(Z_i)}^{p-1} + \abs{(\pi_{s}f)(Z_i)}^{p-1} \Big), \\
    B_s(f) &\lesssim_p u2^{s/2} \biggl( \sum_{i\notin I_s}(\Delta_s f)(Z_i)^2 \Big( \abs{(\pi_{s+1}f)(Z_i)}^{2p-2} + \abs{(\pi_{s}f)(Z_i)}^{2p-2} \Big) \bigg)^{1/2}.
\end{align*}

\smallskip

{\em Step 4: Order statistics.}
Toward applying Lemma \ref{lem: rearrangement}, let $t \coloneqq 2u^2 2^s$ and $k_s \coloneqq t/\ln_+(en/t)$, and choose $I_s \subset [N]$ as a set with cardinality $\lfloor 3k_s \rfloor$ which contains $\lfloor k_s \rfloor$ largest coordinates of $|(\Delta_sf)(Z_i)|$, $\lfloor k_s \rfloor$ largest coordinates of $|(\pi_{s+1}f)(Z_i)|$ and $\lfloor k_s \rfloor$ largest coordinates of $|(\pi_{s}f)(Z_i)|$ (and any other indices in $[N]$ to make the total cardinality $\lfloor 3k_s \rfloor$). Thus, applying Cauchy-Schwarz inequality and recalling that $\pi_{s+1}f,\, \pi_s f \in \mathcal{F}$, we get 
\begin{equation} \label{eq:preliminary_estimate_I}
   A_s(f) \leq 2 \underbrace{\bigg( \sum_{i \le 3k_s} \big( (\Delta_s f)(Z_i)^* \big)^2 \bigg)^{1/2}}_{\eqqcolon A'_s(f)} \underbrace{\sup_{\phi \in \mathcal{F}}\bigg(\sum_{i=1}^N \abs{\phi(Z_i)}^{2(p-1)} \bigg)^{1/2}}_{\eqqcolon \Phi},
\end{equation}
where the star indicates a nonincreasing rearrangement of $N$ numbers (introduced just above Lemma \ref{lem: rearrangement}).  

To bound $A'_s(f)$, we apply Lemma \ref{lem: rearrangement} for $k = k_s$ and for the random variables $X_i = (\Delta_s f)(Z_i) = \ip{\Delta_s v}{Z_i}$, which satisfy $\norm{X_i}_\psitwo \lesssim_K \norm{\Delta_s v}_2$ since $Z_i$ are subgaussian. After rescaling, we get with probability at least $1-2e^{-2u^22^s}$ that
\begin{equation}    \label{eq: first half of As}
    A'_s(f)
    \lesssim_K u2^{s/2}\|\Delta_sv\|_2.
\end{equation}

\smallskip

{\em Step 5: Matrix deviation inequality.}
Recall that $A$ is the random matrix whose rows are $Z_1,\cdots, Z_N$. To bound $\Phi$, recall from the definition of $\mathcal{F}$ that
\begin{equation}    \label{eq: function to matrix}
    \Phi = \sup_{v\in T}\bigg(\sum_{i=1}^N |\langle Z_i,v\rangle|^{2(p-1)}\bigg)^{1/2} = \sup_{v\in T}\norm{Av}_{2(p-1)}^{p-1}.
\end{equation}
By Theorem \ref{thm:matrix_dev} (see Remark \ref{rmk:tail_matrixdev}), with probability at least $1-2e^{-u^2}$, the following event holds:
\begin{equation}
\label{eq:matrix_dev}
    \sup_{v\in T} 
    \abs*{ \norm{Av}_{2(p-1)}-N^{1/2(p-1)} \norm{\ip{Z}{v}}_{L^{2(p-1)}}}
    \lesssim_{K,p} \gamma(T) + u\rad(T).
\end{equation}
Suppose this event occurs. By the following estimate $$\norm{\ip{Z}{v}}_{L^p} \lesssim_{K,p} \norm{\ip{Z}{v}}_{L^2} = \sqrt{v^\tran \E ZZ^\tran v}= \sqrt{v^\tran I_d v}=\norm{v}_2 \le \rad(T),$$ together with the reverse triangle inequality, we have
\begin{equation*}
    \sup_{v\in T} \norm{Av}_{2(p-1)}
    \lesssim_{K,p} \gamma(T) + u\rad(T) + N^{1/2(p-1)}\rad(T).
\end{equation*}
Then, using the numeric inequality $(a+b+c)^{p-1} \leq 3^{p-2} (a^{p-1}+b^{p-1}+c^{p-1})$ which for $p>2$ follows from Jensen's inequality and recalling \eqref{eq: function to matrix}, we conclude that the event
\begin{equation*}    
    \mathcal{E} 
    \coloneqq \left\{ \Phi
    \lesssim_{K,p} \gamma(T)^{p-1} + (\sqrt{N}+u^{p-1})\rad(T)^{p-1} \right\}
\end{equation*}
is likely, i.e, 
\begin{equation} \label{eq: E is likely} 
    \Probevent{\EE} \ge 1-2e^{-u^2}.
\end{equation}
Substituting this and \eqref{eq: first half of As} into \eqref{eq:preliminary_estimate_I}, we conclude that, for any $s \ge 0$ and $f \in \FF$, the following bound holds with probability at least $1-2e^{-2u^22^s}$:
\begin{equation}    \label{eq: bound on As}
    \one_{\EE} A_s(f)
    \lesssim_{K,p} u2^{s/2} \norm{\Delta_s v}_2 \underbrace{\Big( \gamma(T)^{p-1} + (\sqrt{N}+u^{p-1})\rad(T)^{p-1} \Big)}_{\eqqcolon\lambda(T)}.
\end{equation}
For the term $B_s(f)$, it is enough to bound
$$
\tilde{B}_s(f) 
\coloneqq u2^{s/2} \bigg( \sum_{i\notin I_s}(\Delta_s f)(Z_i)^2 \abs{(\pi_s f)(Z_i)}^{2p-2} \bigg)^{1/2},
$$
as the bound with $\pi_{s+1}f$ follows similarly. Applying first the Cauchy-Schwarz inequality and then Lemma \ref{lem: rearrangement} with $q=4$ (where we use again that $\norm{\ip{\Delta_s}{Z_i}}_\psitwo\lesssim_K\norm{\Delta_sv}_2$), we obtain with probability at least $1-2e^{-2u^22^s}$ that
\begin{align*}
    \tilde{B}_s(f) 
    &\le u2^{s/2}\bigg( \sum_{i>k_s}\big( (\Delta_s f)(Z_i)^* \big)^4 \bigg)^{1/4}\bigg(\sum_{i > k_s}\abs*{(\pi_s f)(Z_i)^*}^{4p-4}\bigg)^{1/4} \\
    &u2^{s/2}\bigg(\sum_{i>k_s}(\Delta_s^{\ast}f)^4(Z_i)\bigg)^{1/4}\bigg(\sum_{i>k_s}|\pi_s f(Z_i)^{\ast}|^{4p-4}\bigg)^{1/4}\\
    &\lesssim_{K,p} u2^{s/2} \Big( N \norm{\Delta_s v}_2^{4} \Big)^{1/4} \Big(N \norm{\pi_s v}_2^{4p-4} \Big)^{1/4} \\
    &\lesssim_{K,p} u2^{s/2}\sqrt{N} \|\Delta_sv\|_2 \rad(T)^{p-1}.
\end{align*}
The same bound holds replacing $\pi_{s+1}$ by $\pi_s$. Therefore, with probability at least $1-4e^{-2u^22^s}$, we get
\begin{equation} \label{eq:final_estimate_B_s}
    B_s(f) \lesssim_{K,p} u2^{s/2}\sqrt{N} \norm{\Delta_s v}_2 \rad(T)^{p-1}.
\end{equation}
Combining this with \eqref{eq: bound on As}, we conclude that the following bound holds with probability at least $1-6e^{-2u^22^s}$:
$$
\one_{\EE} \big( A_s(f)+B_s(f) \big)
\lesssim_{K,p} u2^{s/2} \norm{\Delta_s v}_2 \lambda(T).
$$ 
Substituting into \eqref{eq: As Bs}, we obtain that, with probability at least $1-7e^{-2u^22^s}$,
\begin{equation} \label{eq: As+Bs}
    \one_{\EE} \abs[\bigg]{\sum_{i=1}^N \varepsilon_i \big( (\pi_{s+1}f)^p-(\pi_sf)^p \big)(Z_i)}
    \lesssim_{K,p} u2^{s/2} \norm{\Delta_s v}_2 \lambda(T).
\end{equation}

\smallskip

{\em Step 6: Union bound.}
Now take a union bound over all $f \in \FF$, or equivalently over all $v \in T$. By definition of an admissible sequence, the number of different pairs $(\pi_s(f), \pi_{s+1}(f))$ is at most $2^{2^s} \cdot 2^{2^{s+1}} \le 2^{3 \cdot 2^s}$. So, taking the union bound over all these, we conclude that for every $u\ge 2$, with probability at least $1-2^{3 \cdot 2^s} \cdot 7e^{-2u^22^s} \ge 1-7e^{-u^22^s}$, the bound \eqref{eq: As+Bs} holds for all $f \in \FF$ simultaneously.

Next, take a union bound over all $s \ge 0$. We obtain that for $u\ge 2$, with probability at least 
\begin{equation*}
1-\sum_{s\ge 0}7e^{-u^22^s} \ge 1-10e^{-u^2},
\end{equation*}
the bound \eqref{eq: As+Bs} holds for all $f \in \FF$ and $s \ge 0$ simultaneously. Furthermore, recalling from \eqref{eq: E is likely} that the event $\EE$ is likely, we conclude that for $u\ge 2$, with probability at least $1-12e^{-u^2}$, the bound
$$
\abs[\bigg]{\sum_{i=1}^N \varepsilon_i \big( (\pi_{s+1}f)^p-(\pi_sf)^p \big)(Z_i)}
\lesssim_{K,p} u2^{s/2} \norm{\Delta_s v}_2 \lambda(T)
$$
holds for all $f \in \FF$ and $s \ge 0$ simultaneously. 

Sum up these inequalities over all $s \ge 0$ and take supremum over $f \in \FF$ (or equivalently, over $v \in T$) on both sides. For $u\ge 2$ we get that with probability at least $1-12e^{-u^2}$:
\begin{align*}
    \sup_{f \in \FF} &\abs[\bigg]{\sum_{s \ge 0} \sum_{i=1}^N \varepsilon_i \big( (\pi_{s+1}f)^p-(\pi_sf)^p \big)(Z_i)} \\
    &\lesssim_{K,p} u \lambda(T) \sup_{v \in T} \sum_{s \ge 0} 2^{s/2} \norm{\Delta_s v}_2 \\
    &\lesssim u \lambda(T) \gamma(T)        
        \quad \text{(by the choice of admissible sequence in \eqref{eq:Talagrand_MM})} \\
    &\lesssim u^p \big( \gamma(T)^p + \sqrt{N} \gamma(T) \rad(T)^{p-1} \big)
        \quad \text{(by definition of $\lambda(T)$ in \eqref{eq: bound on As})}.
\end{align*}
Integrating the tail, we conclude that 
\begin{equation} \label{eq: sup term}
    \E \sup_{f \in \FF} \abs[\bigg]{\sum_{s \ge 0} \sum_{i=1}^N \varepsilon_i \big( (\pi_{s+1}f)^p-(\pi_sf)^p \big)(Z_i)}
    \lesssim_{K,p} \gamma(T)^p + \sqrt{N} \gamma(T) \rad(T)^{p-1}.
\end{equation}

\smallskip

{\em Step 7: The first term of chaining.}
We successfully bounded the first term in \eqref{eq: chained}. The second term
is much easier, since there is no supremum:
\begin{align}
    \E \abs[\bigg]{\sum_{i=1}^N\varepsilon_i(\pi_{0}f)^p(Z_i)}
    &= \E_Z\E_{\varepsilon} \abs[\bigg]{\sum_{i=1}^N \varepsilon_i \langle Z_i,\pi_{0}v\rangle^p}\\
    &\le \E \bigg( \sum_{i=1}^N\langle Z_i,\pi_0v\rangle^{2p} \bigg)^{1/2} \nonumber \quad \text{(Jensen's inequality)}\\
    &\lesssim_{K,p} \sqrt{N}\rad(T)^{p}
    \lesssim \sqrt{N}\gamma(T)\rad(T)^{p-1}. \label{eq: fixed term}
\end{align}
The first inequality follows by first taking expectation with respect to the random signs $(\varepsilon_i)$ and then with respect to $(Z_i)$. To get the second inequality, follow the steps we made to bound $\Phi$ in \eqref{eq: E is likely}, using matrix deviation inequality \eqref{eq:matrix_dev} with $2p$ instead of $2(p-1)$. The third inequality follows since we always have $\rad(T) \lesssim \gamma(T)$.

Substitute \eqref{eq: sup term} and \eqref{eq: fixed term} into \eqref{eq: chained} and then into \eqref{eq: symmetrization}, and 
Theorem \ref{thm:main} is proved.
\end{proof}

\section{Remarks}
The proof does not rely on $p$ being an integer. The same bound (up to an absolute constant) holds for the process
\begin{equation*}
\sup_{v\in T}\left|\frac{1}{N}\sum_{i=1}^N |\langle Z_i,v\rangle|^p - \E|\langle Z,v\rangle|^p \right|,
\end{equation*}
where $p\ge 2$ may not be an integer. We assumed that $p$ is an integer in the case of simple tensors just because the function $x^p$ might not be well-defined for non-integer $p$. 

Similarly as in the proof of Al-Ghattas, Chen and Sanz-Alonso, our proof extends to any star-shaped class of functions $\mathcal{F}$ such that for any $f\in \mathcal{F}$, it holds that $\|f(Z)-g(Z)\|_{\psi_2}\lesssim \|f(Z)-g(Z)\|_{L^2}$.


\begin{appendices}

\section{Proof of Lemma \ref{lem: rearrangement}}\label{a: order statistics}

{\em Case 1: $t \le n$.} In this regime, we have
\begin{equation}    \label{eq: t vs k}
    k \le n
    \quad \text{and} \quad
    t \asymp k \ln(en/k).
\end{equation}
The first bound follows if we replace $t$ with $n$ in the definition of $k$. The second bound follows from choosing $k = t/\ln(en/t)$. 

To prove the first bound in \eqref{eq: order statistics}, note that if $\sum_{i \le k} (X_i^*)^2 > a$, then there exists a subset $I \subset [n]$ with $\abs{I} \le 3k$ for which $\sum_{i \in I} X_i^2 > a$. There are at most $(en/k)^{3k}$ subsets of cardinality bounded by $3k$, so the union bound implies that there exists an absolute constant $C>0$ for which it holds that
\begin{equation} \label{eq: union bound first coeffs}
    \Prob[\Bigg]{\sum_{i \le k} (X_i^*)^2 > C(k+t)}
    \le \Big(\frac{en}{k} \Big)^{3k} \max_{\abs{I}=k} \Prob[\Bigg]{\sum_{i \in I} X_i^2 > C(k+t)},
\end{equation}
and that $\E X_i^2 \le C \norm{X_i}_\psitwo^2 \le C$. So, 
$$
\Prob[\Bigg]{\sum_{i \in I} X_i^2 > C(k+t)}
= \Prob[\Bigg]{\sum_{i \in I} (X_i^2 - \E X_i^2) > Ct}.
$$
The random variables $X_i^2 - \E X_i^2$ are independent, mean zero, and subexponential: $\norm{X_i^2 - \E X_i^2}_\psione \lesssim \norm{X_i^2}_\psione = \norm{X_i}_\psitwo^2 \le 1$, so Bernstein's inequality \cite[Theorem~2.9.1]{vershynin2018high} implies that there exists a constant $c>0$ such that
$$
\Prob[\Bigg]{\sum_{i \in I} (X_i^2 - \E X_i^2) > Ct}\le 2 \exp \Big( -c \min((Ct)^2/k, Ct) \Big)
\le 2 \exp(-cCt),
$$
where we used that $t \gtrsim k$ from \eqref{eq: t vs k} and chose $C$ sufficiently large. Hence the probability in \eqref{eq: union bound first coeffs} is bounded by 
$$
2(en/k)^{3k} \exp(-cCt)
\le 2\exp(-t),
$$
where we again used \eqref{eq: t vs k} and chose $C$ sufficiently large. Thus, with probability at least $1-2e^{-t}$, we have 
$$
\sum_{i \le k} (X_i^*)^2 
\lesssim k+t
\lesssim t,
$$
where the last bound follows from \eqref{eq: t vs k}. The first bound in \eqref{eq: order statistics} is proved. 

To prove the second bound in \eqref{eq: order statistics}, choose a sufficiently large absolute constant $C>0$, fix any $i \in [n]$ and argue like above to get 
\begin{align*}
    \Prob[\Bigg]{X_i^* > C\sqrt{\ln\frac{en}{i}}}
    &\le \binom{n}{i} \max_{i \in [n]} \Big( \Prob[\Big]{X_i > C\sqrt{\ln\frac{en}{i}}} \Big)^i \\
    &\le \Big(\frac{en}{i}\Big)^i \Big(\exp \big(-2\ln\frac{en}{i}\big)\Big)^i
    = \Big(\frac{en}{i}\Big)^{-i}.
\end{align*}
Thus
$$
\Prob[\Bigg]{\exists i>k:\; X_i^* > C\sqrt{\ln\frac{en}{i}}}
\le \sum_{i>k} \Big(\frac{en}{i}\Big)^{-i}
\lesssim \Big(\frac{en}{k}\Big)^{-k}
\le e^{-t}.
$$
Therefore, with probability at least $1-e^{-t}$,  we have 
$$
\sum_{i>k} (X_i^*)^q
\lesssim \sum_{i=1}^n \Big( \ln\frac{en}{i} \Big)^{q/2}
\lesssim n,
$$
proving the second bound in \eqref{eq: order statistics}.

\medskip

{\em Case 2: $t>n$.} In this regime, we have $k>n$, which can be seen by replacing $t$ with $n$ in the definition of $k$. So, only the first bound in \eqref{eq: order statistics} needs to be shown. We have
$$
\Prob[\Bigg]{\sum_{i \le 3k} (X_i^*)^2 > C(n+t)}
= \Prob[\Bigg]{\sum_{i=1}^n X_i^2 > C(n+t)}
\le 2\exp(-t),
$$
where the last bound follows in a way similar to Case 1 (but with $n$ instead of $k$). Thus, with probability at least $1-2e^{-t}$, we have 
$$
\sum_{i \le 3k} (X_i^*)^2 
\lesssim n+t
\lesssim t.
$$
The first bound in \eqref{eq: order statistics} is proved, finishing the proof of Lemma \ref{lem: rearrangement}.




\end{appendices}

\section*{Acknowledgments}
PA and RV are supported by NSF and Simons Research Collaborations on the Mathematical and Scientific Foundations of Deep Learning. RV is also supported by NSF Grant DMS 1954233.

\end{document}